\documentclass[11pt,leqno]{amsart}

\usepackage{amssymb}
\usepackage{amsmath}
\usepackage{amsxtra}
\usepackage{amscd}
\usepackage{amsfonts}
\usepackage[utf8]{inputenc}
\usepackage{hyperref}
\usepackage[all]{xy}
\usepackage{tikz}
\usetikzlibrary{arrows, decorations.markings}
\usepackage{float}
\usepackage[vcentermath]{youngtab}

\newtheorem{theorem}{Theorem}
\newtheorem{corollary}{Corollary}

\newtheorem{proposition}{Proposition}

\numberwithin{equation}{section}

\makeatletter
\let\@wraptoccontribs\wraptoccontribs
\makeatother

\begin{document} 

\title{A Survey on Temperley-Lieb-type quotients from the Yokonuma-Hecke algebras}

\author{D. Goundaroulis}
\address{Center Int\'egratif de G\'enomique,
UNIL,
Batiment G\'enopode, CH-1015 Lausanne, Switzerland.}
\email{dimoklis.gkountaroulis@unil.ch}
%\urladdr{users.ntua.gr/dgound}



\date{}

\begin{abstract}
In this survey we collect all results regarding the construction of the Framization of the Temperley-Lieb algebra of type $A$ as a quotient algebra of the Yokonuma-Hecke algebra of type $A$. More precisely, we present all three possible quotient algebras the emerged during this construction and we discuss their dimension, linear bases, representation theory and the necessary and sufficient conditions for the unique Markov trace of the Yokonuma-Hecke algebra to factor through  to each one of them. Further, we present the link invariants that are derived from each quotient algebra and we point out which quotient algebra provides the most natural definition for a framization of the Temperley-Lieb algebra. From the Framization of the Temperley-Lieb algebra we obtain new one-variable invariants for oriented classical links that, when compared to the Jones polynomial, they are not topologically equivalent since they distinguish more pair of non isotopic oriented links.  Finally, we discuss the generalization of the newly obtained  invariants to a new two-variable invariant for oriented classical links that is stronger than the Jones polynomial.
\end {abstract}
\maketitle

\section{Introduction}
The Yokonuma-Hecke algebra was first introduced in the 60's by Yokonuma as a generalization of the Iwahori-Hecke algebra in the context of Chevalley groups \cite{yo}. In recent years, Juyumaya  simplified the natural description by giving a presentation in terms of generators and relations \cite{jukan,ju,jusur}. A detailed overview of Juyumaya's approach can be found in \cite[Preliminaries]{marin}. In this context, the Yokonuma-Hecke algebra of type $A$ can be considered as a quotient of the framed braid group algebra over a two-sided ideal that is generated by a quadratic relation that involves certain weighted idempotent elements. 

Throughout the past ten years the theory of Yokonuma-Hecke algebras received a significant amount of attention, mainly due to the concept of framization of knot algebras, a concept that was introduced by Juyumaya and Lambropoulou \cite{jula}.  A knot algebra is an algebra that is involved in the construction of invariants of classical links via braid group representations \cite{jo}. To be more precise, a knot algebra ${\rm A}$ is a triplet $({\rm A}, \pi , \tau ) $, where $\pi$ is an appropriate representation of the braid group in ${\rm A}$ and $\tau$ is a Markov trace function defined on ${\rm A}$. The Iwahori-Hecke algebra and the Temperley-Lieb algebra are the most known examples of knot algebras. On the other hand, the framization consists in an extension of  a knot algebra via the  addition of  framing generators which gives rise to a new algebra that is related to framed braids and framed knots.  The Yokonuma-Hecke algebra, ${\rm Y}_{d,n}(u)$ is the basic example of this concept and it can be regarded as a  framization of the Iwahori-Hecke algebra, ${\rm H}_n(u)$ \cite{jula2, jula}. With this in mind, Juyumaya and Lambropoulou proposed framizations of several knot algebras \cite{jula5, jula6} from which  isotopy invariants for framed, classical and singular links  were derived \cite{jula2,jula4,jula3}.

The breakthrough in this theory came while comparing the invariants for classical oriented knots and links from the Yokonuma-Hecke algebras to the Homflypt polynomial. In \cite{ChJuKaLa}, the use of a different presentation for ${\rm Y}_{d,n}$ with parameter $q$ instead of $u$ and a different quadratic relation led to the proof that the derived two-variable invariants $\Theta_d$ are {\it not topologically equivalent to the Homflypt polynomial on links} while they are topologically equivalent to the Homflypt on knots. Furthermore, in the same work it was shown that the invariants $\Theta_d$ distinguish more pairs of non-isotopic oriented links than the Homflypt polynomial. Moreover, it was shown that the invariants can be generalized to a 3-variable invariant $\Theta$ for oriented classical links that can be completely defined via the skein relation of the Homflypt polynomial on crossings involving different components of the link and a set of initial conditions \cite{ChJuKaLa, kaula}. The invariant $\Theta$ distinguishes the same pairs of Homflypt-equivalent links as $\Theta_d$, it is not topologically equivalent to the Homflypt or the Kauffman polynomials and, thus, {\it it is stronger than the Homflypt polynomial on links}.

One of the open problems in the concept of framization of knot algebras was the determination of a framization of the  Temperley-Lieb algebra. If one considers the classical Temperley-Lieb algebra as was introduced by Jones \cite{jo} that is, as a quotient of the Iwahori-Hecke algebra, it is immediately evident that desired framization will emerge as an appropriate quotient of the Yokonuma Hecke algebra. Contrary to the classical case such a candidate algebra is not unique. The study of these quotient algebras has been the topic of the author's PhD thesis \cite{go} which led to a series of results regarding their topological \cite{gojukola, gojukola2, gola, gola1} as well as their algebraic properties \cite{ChPou, ChPou2}.  There are three potential candidates that can qualify as the framization of the Temperley-Lieb algebra: the {\it Yokonuma-Temperley-Lieb algebra} ${\rm YTL}_{d,n}(u)$,  the {\it Complex Reflection Temperley-Lieb algebra} ${\rm CTL}_{d,n}(u)$ and the so-called {\it Framization of the Temperley-Lieb algebra} ${\rm FTL}_{d,n}(u)$.  The algebra ${\rm YTL}_{d,n}(u)$ is too restricted and, as a consequence,  the invariants for classical links from the algebra ${\rm YTL}_{d,n}(u)$ just recover the Jones polynomial \cite{gojukola}. On the other hand, the algebra ${\rm CTL}_{d,n}(u)$ is too large for our topological purposes and the derived link invariants coincide either with those from ${\rm Y}_{d,n}(u)$ or with those from ${\rm FTL}_{d,n}(u)$ \cite{gojukola2}. Unfortunately, these two quotient algebras do not fit the topological purposes of deriving new invariants for (framed) knots and links and, thus, they do not qualify as potential framizations of the Temperley-Lieb algebra. The third quotient algebra of ${\rm Y}_{d,n}(u)$, the Framization of the Temperley-Lieb algebra, ${\rm FTL}_{d,n}(u)$, lies between ${\rm YTL}_{d,n}(u)$ and ${\rm CTL}_{d,n}(u)$ and, as it will be made clear in Section~\ref{threequot}, it turns out to be the right one \cite{gojukola2}. The invariants $\theta_d$ for classical links from the algebras ${\rm FTL}_{d,n}$ adapted to a presentation with parameter $q$ instead of $u$ of ${\rm Y}_{d,n}(q)$, are proven to be {\it not topologically equivalent to the Jones polynomial on links} while they are topological equivalent to the Jones polynomial on knots \cite{gojukola2}. Finally, in analogy to the invariants $\Theta_d$, the invariants $\theta_d$ can be generalized to a new two-variable invariant of oriented classical links $\theta$ that is {\it stronger than the Jones polynomial} \cite{gola1}.

The outline of the paper is as follows:
In Section~\ref{prelim} we introduce the necessary notations and we give a brief overview of all the required definitions and results such as: the Temperley-Lieb algebra, the Yokonuma-Hecke algebra, the ${\rm E}$-system and the derived two-variable invariants for oriented framed and classical knots and links.  In Section~\ref{threequot} we discuss three quotients of the Yokonuma-Hecke algebra as possible candidates for the framization of the Temperley-Lieb algebra. Moreover, we give all algebraic (linear basis, dimension, representation theory) as well as all topological (Markov trace, link invariants) results in the literature regarding each one of these quotient algebras. In Section~\ref{compgen} we describe how the invariants $\Theta_d$ and $\theta_d$ compare to the Homflypt and Jones polynomials respectively. Finally, we discuss how the invariants $\Theta_d$ generalize to a new three-variable invariant for oriented classical links as well as the analogous generalization of the invariants $\theta_d$ to a new two-variable invariant for oriented classical links and we also describe closed combinatorial formulas for each one of the generalizations.

\section{Preliminaries}\label{prelim}
In this section we will establish our notation and we will present the basic notions that will be used in the following sections. 
\subsection{{\it Notations}} We start by fixing two positive integers, $d$ and $n$. Every algebra considered in this paper is an associative unital algebra over the field $\mathbb{C}(u)$, where $u$ is an indeterminate.   The {\it framed braid group} on $n$ strands is defined as the semi-direct product of Artin's braid group $B_n$ with $n$ copies of $\mathbb{Z}$, namely: ${\mathcal F}_{n} = \mathbb{Z}^n \rtimes  B_n$, where the action of the braid group $B_n$ on $\mathbb{Z}^n$ is given by the permutation induced by a braid on the indices $\sigma_it_j=t_{s_i(j)}\sigma_i$.  By considering framings modulo $d$, the {\it modular framed braid group}, $\mathcal{F}_{d,n} = (\mathbb{Z}/d\mathbb{Z})^n \rtimes B_n$, is defined. Due to the above action a word $w$ in ${\mathcal F}_{n}$ (resp. $\mathcal{F}_{d,n}$) has the {\it splitting property}, that is, it splits into the  {\it framing} part and the {\it braiding}  part
$
w = t_1^{a_1}\ldots t_n^{a_n} \, \sigma
$
where $\sigma \in B_n$ and $a_i \in \mathbb{Z}$ (resp. $\mathbb{Z}/d\mathbb{Z}$). 

Finally, {\it a partition of} $n$, $\lambda = (\lambda_1, \ldots , \lambda_k)$, is a  family of positive integers such that $\lambda_1 \geq \lambda_2 \geq \ldots \geq \lambda_k \geq 1$ and $|\lambda|=\lambda_1 + \ldots+ \lambda_k = n$. We identify every partition with its Young diagram, that is a left-justified array of $k$ rows such that the $j$-th row contains $\lambda_j$ nodes for all $j = 1, \ldots, k$. A $d$-partition $\lambda$, or a Young $d$-diagram, of size $n$ is a $d$-tuple of partitions such that the total number of nodes in the associated Young diagrams is equal to $n$. That is, we have $\lambda = (\lambda^{(1)},\ldots ,\lambda^{(d)})$  with $\lambda^{(1)},\ldots ,\lambda^{(d)},$  usual partitions such that $|\lambda^{(1)}|+ \ldots + |\lambda^{(d)}| = n$.

\subsection{{\it The Temperley-Lieb algebra}} For $n \geq 3$, the Temperley-Lieb algebra, ${\rm TL}_n(u)$, is the $\mathbb{C}(u)$-algebra that is generated by the elements $h_1, \ldots , h_{n-1}$ which are subject to the following relations:
 \begin{align*}
h_i h_j &= h_j h_i \quad \text{for all} \quad |i-j| >1 \\
h_i h_j h_i &= h_j h_i h_j \quad \text{for all} \quad \vert i - j \vert =1 \\
h_i^2 &=  u + (u-1) h_i \\
h_{i,i+1} &= 0,
\end{align*}
where $h_{i,j}:=1+ h_i + h_j +h_i h_j + h_j h_i + h_i h_j h_i$. Notice that the first three relations are the defining relations of the Iwahori-Hecke algebra, ${\rm H}_n(u)$, which is defined as the quotient of the algebra $\mathbb{C}(u)B_n$ over the two-sided ideal that is generated by the quadratic relations mentioned above. Thus, with this presentation, the algebra ${\rm TL}_n(u)$ can be considered as the quotient of ${\rm H}_n(u)$ over the two-sided ideal that is generated by the elements $h_{i,i+1} \in {\rm H}_n(u)$.  It is not difficult to see that the defining ideal of ${\rm TL}_{n}(u)$ is principal and that is generated by the element $h_{1,2}$. 

The algebra ${\rm H}_n(u)$ supports a unique Markov trace, the Ocneanu trace $\tau$ with parameter $\zeta$ \cite[Theorem~5.1]{jo}. By normalizing and rescaling $\tau$ according to the braid equivalence, one obtains the {\it Homlypt polynomial} \cite[Proposition~6.2]{jo}, \cite{homfly, pt}. Further, the trace $\tau$ factors through to the quotient algebra ${\rm TL}_n(u)$. The necessary and sufficient conditions for the factoring of $\tau$ provide a specialization for the trace parameter $\zeta$ which, in turn, gives rise to the {\it Jones polynomial} \cite{jo}:
\[
V(u)(\widehat{\alpha}) = \left(- \frac{1+u}{\sqrt{u}} \right)^{n-1} \left(\sqrt{u}\right)^{\varepsilon(\alpha)} {\rm \tau}(\pi(\alpha)),
\]
where: $\alpha \, \in \, \cup_{\infty} B_{n}$, $\pi$ is the natural epimorphism of $\mathbb{C}(u)B_n$ on ${\rm TL}_n(u)$ that sends the braid generator $\sigma_i$ to $h_i$ and $\varepsilon(\alpha)$ is the algebraic sum of the exponents of the $\sigma_i$'s in $\alpha$.
\subsection{{\it The Yokonuma-Hecke algebra}}\label{yhsec} The {\it Yokonuma-Hecke algebra} ${\rm Y}_{d,n}(u)$ \cite{yo} is defined as the quotient of the group algebra
$\mathbb{C}(u) {\mathcal F}_{d,n}$ over the two-sided ideal  generated by the elements:
\[
\sigma_i^2  -  1 - (u-1)e_i - (u-1)e_i \sigma_i \quad \mbox{for all } i,
\]
where $e_{i} := \frac{1}{d} \sum_{s=0}^{d-1}t_i^s t_{i+1}^{d-s}$, for  $i=1,\ldots , n-1$. The elements $e_i$ in ${\rm Y}_{d,n}(q)$ are idempotents \cite{ju}. The generators of the ideal give rise to the following quadratic relations in ${\rm Y}_{d,n}(q)$: 
\begin{equation}\label{yqeq}
g_i^2 = 1 + (u-1) e_i + (u-1)e_i g_i,
\end{equation}
where $g_i$ corresponds to $\sigma_i$. Moreover, \eqref{yqeq} implies that the elements $g_i$ are invertible with $g_i^{-1} = g_i - (u^{-1}- 1) e_i +(u^{-1}-1)e_ig_i$, $ 1 \leq i \leq n-1$. The $t_i$'s are called the {\it framing generators}, while the $g_i$'s are called the {\it braiding generators} of ${\rm Y}_{d,n}(q)$. By its construction, the Yokonuma-Hecke algebra is considered as {\it the framization of the Iwahori-Hecke algebra}. Regarding its algebraic properties, the algebra ${\rm Y}_{d,n}(u)$ has the following  standard linear basis  \cite{ju}:
\[
\{t_1^{a_1}\ldots t_n^{a_n} w \, | \,a_i \in \mathbb{Z}/d\mathbb{Z},\,  w \in \mathcal{B}_{{\rm H}_n} \},
\]
where  $\mathcal{B}_{{\rm H}_n}$ is the standard basis of ${\rm H}_n(u)$. A simple counting argument implies that the dimension of the algebra ${\rm Y}_{d,n}(u)$ is equal to $n!\ d^n $. Further, the irreducible representations of $Y_{d,n}(u)$ over $\mathbb{C}(u)$, are parametrised by the $d$-partitions of $n$ \cite[Theorem~1]{ChPoulain}.

One of the most important results regarding  the Yokonuma-Hecke algebra lies in \cite{ju} where Juyumaya showed that  ${\rm Y}_{d,n}(u)$ supports the following unique linear Markov trace function: 
\[
{\rm tr}_d:  \cup_{n=1}^{\infty}{\rm Y}_{d,n}(u) \longrightarrow   \mathbb{C}(u)[z, x_1, \ldots, x_{d-1}],
\]
 where $z$, $x_1$, $\ldots, x_{d-1}$ are indeterminates. The trace ${\rm tr}_d$ can be defined inductively on $n$ by the following rules \cite[Theorem~12]{ju}:
\[
\begin{array}{rcll}
{\rm tr}_d(ab) & = & {\rm tr}_d(ba)  \qquad &  \\
{\rm tr}_d(1) & = & 1 & \\
{\rm tr}_d(ag_n) & = & z\, {\rm tr}_d(a) \qquad & \\
{\rm tr}_d(at_{n+1}^s) & = & x_s {\rm tr}_d(a)\qquad  & (  s = 1, \ldots , d-1) ,
\end{array}
\]
where  $a,b \in {\rm Y}_{d,n}(u)$. Using the rules of ${\rm tr}_d$ and setting  $x_0:=1$, one deduces that ${\rm tr}_d(e_i)$ takes the same value for all $i$, indeed: $E := {\rm tr}_d(e_i)= \frac{1}{d}\sum_{s=0}^{d-1}x_{s}x_{d-s}.$ 

In order to define framed and classical link invariants via the trace ${\rm tr}_d$, one should re-scale ${\rm tr}_d$ according to the framed braid equivalence \cite{ks}. Unfortunately, the trace ${\rm tr}_d$ is the only known trace that does not re-scale directly \cite{jula}. The {\it ${\rm E}$-system} is the following system of non-linear equations
\[
\sum_{s=0}^{d-1}x_{m+s}x_{d-s} =   x_{m}\sum_{s=0}^{d-1}x_{s}x_{d-s} \qquad (1\leq m \leq d-1),
\]
 that was introduced in order to find the  necessary and sufficient conditions that needed to be applied on the parameters $x_i$  of  $\rm tr$ so that the re-scaling of ${\rm tr}_d$ would be possible \cite{jula}. We say that the $(d-1)$-tuple  of complex numbers $({\rm x}_1, \ldots , {\rm x}_{d-1})$ satisfies the {\it ${\rm E}$-condition} if ${\rm x}_1,\ldots , {\rm x}_{d-1}$ are solutions of the ${\rm E}$-system. The full set of solutions of the ${\rm E}$-system is given by Paul G\'{e}rardin \cite[Appendix]{jula} using tools of harmonic analysis on finite groups. More precisely, he interpreted the solution 
{$({\rm x}_1,\ldots,{\rm x}_d)$} 
of the ${\rm E}$-system, as the complex function $x : \mathbb{Z}/d\mathbb{Z} \rightarrow \mathbb{C}$ that sends $k \mapsto {\rm x}_k$, $k\neq 0$ and $0 \mapsto 1$. Let now $\chi_m$ be the character of the group $\mathbb{Z}/d\mathbb{Z}$ and let $\mathbf{i}_m:=\sum_{s=0}^{d-1} \chi_m( s) t^s$, for $ m\in {\mathbb Z}/d{\mathbb Z} \in \mathbb{C}[\mathbb{Z}/d\mathbb{Z}]$. We then have that the solutions of the ${\rm E}$-system are of the following form:
\[
{\rm x}_s = \frac{1}{|D|} \sum_{m \in D} \mathbf{i}_m (s), \quad 1 \leq s \leq d-1,
\]
where  $D$ is a non-empty subset of $\mathbb{Z}/d\mathbb{Z}$. Hence, the solutions of the ${\rm E}$-system are parametrized by the non-empty subsets of $\mathbb{Z}/d\mathbb{Z}$. Two obvious solutions of the ${\rm E}$-system are: the all-zero solution, that is  $x_i=0$, for all $i$, and when the $x_i$'s are specialized to the $d$-th roots of unity. For the rest of the paper we fix $X_D = \{ {\rm x}_1 , \ldots , {\rm x}_{d-1} \}$ to be a solution of the ${\rm E}$-system parametrized by the non-empty subset $D$ of $\mathbb{Z}/d\mathbb{Z}$. If we specialize the trace parameters $x_i$ of ${\rm tr}_d$ to the values ${\rm x}_i$  we obtain the {\it specialized trace} ${\rm tr}_{d,D}$ with parameter $z$ \cite{ChLa, ChJuKaLa}.  

By normalizing and re-scaling the specialized trace ${\rm tr}_{d,D}$, invariants for {\it framed links} are obtained \cite{jula}:
\begin{equation}\label{gammainv}
\Gamma_{d,D}(w,u)(\widehat{\alpha}) = \left(- \frac{(1 - w u)|D|}{\sqrt{w} (1-u) } \right)^{n-1} \left(\sqrt{w}\right)^{\varepsilon(\alpha)} {\rm tr}_{d,D}(\gamma(\alpha)),
\end{equation}
where: $w = \frac{z + (1-u)}{uz |D|}$ is the re-scaling factor, $\gamma$ is the natural epimorphism of the framed braid group algebra $ \mathbb{C}(u)\mathcal{F}_n$ on the algebra ${\rm Y}_{d,n}(u)$, and $\alpha \in \cup_{\infty} \mathcal{F}_{n}$.   Further, by restricting the invariants $\Gamma_{d,D}(w,u)$ to {\it classical links}, seen as framed links with all framings zero, in \cite{jula4} invariants of classical oriented links $\Delta_{d,D}(w,u)$ are obtained. In  \cite{ChLa} it was proved that for generic values of the parameters $u,z$ the invariants $\Delta_{d,D}(w,u)$ do not coincide with the Homflypt polynomial except in the trivial cases $u=1$ and $E_D=1$. 

\section{The three possible candidates}\label{threequot}
 In this section we will present all results in the literature regarding the three possible quotient algebras that can be considered as candidates for the framization of the Temperley-Lieb algebra.  In what follows, we will give the definitions and dimensions for each quotient algebra, describe their linear bases and representation theory and discuss the necessary and sufficient conditions so that the trace ${\rm tr}_d$ passes to each one of the quotient algebras. Finally, we will present the invariants for framed and classical links that are derived from each algebra.

 \subsection{{\it Motivation behind the construction}}  Following the construction of the classical Temperley-Lieb algebra we would like to introduce an analogue of ${\rm TL}_n(u)$ in the context of framed knot algebras. Namely, to define a quotient of ${\rm Y}_{d,n}(u)$ over a two-sided ideal that is constructed from an appropriately chosen subgroup of the underlying group $C_{d,n}:=\left(\mathbb{Z}/d\mathbb{Z}\right)^n \rtimes S_n$ of ${\rm Y}_{d,n}(u)$. At this point two such subgroups emerge naturally. The first possibility is to consider the subgroups  $\langle s_i , s_{i+1} \rangle $ of $C_{d,n}$ that are also related to the defining ideal of ${\rm TL}_n(u)$. The second possibility is to let the framing generators $t_i$ be involved in the generating set of such a subgroup and consider the following subgroup of $C_{d,n}$:
\[
C^{i}_{d,n} := \langle t_i, t_{i+1}, t_{i+2}\rangle \rtimes \langle  s_i, s_{i+1}\rangle\quad \text{for all } i. 
\] 
Therefore we can define at least two types  of algebras which could be considered  as analogues 
of the Temperley-Lieb algebras in the context of knot algebras with framing. The algebra that corresponds to the first possibility is the {\it Yokonuma-Temperley-Lieb algebra}, denoted by ${\rm YTL}_{d,n}(u)$, while the second is the {\it Complex Reflection Temperley-Lieb algebra}, ${\rm CTL}_{d,n}(u)$.

As mentioned in \eqref{gammainv}, new two-variable invariants for oriented framed knots and links are defined through the trace ${\rm tr}_d$ on the Yokonuma-Hecke algebra by imposing the ${\rm E}$-system on the parameters $x_1, \ldots , x_{d-1}$ \cite{jula}. Hence, we expect that the framization of the Temperley-Lieb algebra will allow us to define one-variable specializations of the invariants derived from ${\rm Y}_{d,n}(u)$. Unfortunately, both quotients above are not satisfactory for this purpose. In the case of ${\rm YTL}_{d,n}(u)$, very strong conditions on the trace parameters must be applied in order for ${\rm tr}_d$ to pass through to the quotient algebra. Namely, the trace parameters $x_i$ must be $d^{th}$ roots of unity, giving rise to obvious, special solutions of the ${\rm E}$-system, which imply topologically loss of the framing information. However, the original Jones polynomial can be recovered from this quotient algebra. In the case of  ${\rm CTL}_{d,n}(u)$, the quotient algebra is large enough so that the necessary and sufficient conditions such that ${\rm tr}_d$ passes to ${\rm CTL}_{d,n}(u)$ are, contrary to the case of ${\rm YTL}_{d,n}(u)$, too relaxed, especially on the trace parameters $x_i$. So, in order to obtain knot invariants we would still need to impose the ${\rm E}$-system on the trace parameters $x_1 , \ldots , x_{d-1}$ as in the case of ${\rm Y}_{d,n}(u)$. 
\smallbreak
The discussion above indicated that the desired framization of the Temperley-Lieb algebra, for our topological purposes, could be an intermediate algebra between the quotient algebras ${\rm YTL}_{d,n}(u)$ and ${\rm CTL}_{d,n}(u)$. One may achieve this, by using for the defining ideal an intermediate subgroup that lies between  $\langle s_i, s_{i+1} \rangle $ and $C^i_{d,n}$. More precisely, we define this  framization  as a quotient of the Yokonuma-Hecke algebra over an ideal that is constructed from the following subgroup of $C_{d,n}$:
\[
H^{i}_{d,n} := \langle t_i t_{i+1}^{-1}, t_{i+1}t_{i+2}^{-1} \rangle \rtimes \langle s_i, s_{i+1} \rangle \quad \text{for all } i.
\]
Thus, one obtains the so-called {\it Framization of the Temperley-Lieb algebra}, ${\rm FTL}_{d,n}(u)$. The  relation between the three quotient algebras is given by the following commutative diagram of epimorphisms \cite[Proposition~3]{gojukola2}:
\[
\xymatrix{ 
&{\rm Y}_{d,n}(u) \ar[d] \ar[r] & {\rm CTL}_{d,n}(u)\ar[d] \ar[r] & {\rm FTL}_{d,n}(u) \ar[dl] \ar[r] & {\rm YTL}_{d,n}(u) \ar@/^/[dll]\\
& {\rm H}_n(u) \ar[r] & {\rm TL}_n (u) & &}
\]
The Yokonuma-Temperley-Lieb algebra and its derived invariants were introduced and studied in \cite{gojukola}, while its representation theory was studied in \cite{ChPou}. The algebras ${\rm FTL}_{d,n}(u)$, ${\rm CTL}_{d,n}(u)$ and their corresponding invariants were introduced in \cite{gojukola2} and were further studied in \cite{ChPou2, gola, gola1}. 

\subsection{{\it The Yokonuma-Temperley-Lieb algebra}}   For $n \geq 3$, the {\it Yokonuma-Temperley-Lieb algebra}, denoted by ${\rm YTL}_{d,n}(u)$, is defined as the quotient of ${\rm Y}_{d,n}(u)$ over the two-sided ideal that is generated by the elements: 
\begin{equation}\label{g12exp}
g_{i,i+1}:= 1+ g_i + g_{i+1} +g_i g_{i+1} + g_{i+1} g_i + g_i g_{i+1} g_i.
\end{equation}

It is a straightforward computation to show that the defining ideal of ${\rm YTL}_{d,n}(u)$ is principal and is generated by the element $g_{1,2}$ \cite[Lemma~4]{gojukola}. Thus, the algebra ${\rm YTL}_{d,n}(u)$ can be considered as the $\mathbb{C}(u)$-algebra that is generated by the elements $t_1, \ldots , t_n , g_1 \ldots ,g_{n-1}$ that are subject to the defining relations of ${\rm Y}_{d,n}(u)$ and the relation $g_{1,2}=0$ \cite[Corollary~1]{gojukola}. Note also that for $d=1$ the algebra ${\rm YTL}_{1,n}(u)$ coincides with ${\rm TL}_n(u)$.

 Every word in the algebra ${\rm YTL}_{d,n}(u)$ inherits the splitting property from ${\rm Y}_{d,n}(u)$. For each fixed element in the braiding part, a set of linear dependency relations among the framing parts can be described which, in turn, lead to the extraction of a linear basis for ${\rm YTL}_{d,n}(u)$ \cite{ChPou}.  Using this technique, Chlouveraki and Pouchin proved in \cite{ChPou} that, for $n \geq 3$, the following set is a linear basis for ${\rm YTL}_{d,n}(u)$:
\[
S_{d,n} = \left \{ t_1^{r_1} \ldots  t_n^{r_n} w \ \vert \ w \in \mathcal{B}_{{\rm TL}} , \ (r_1, \ldots , r_n ) \in \mathcal{E}_{d,n}(w) \right \},
\]
where $\mathcal{B}_{{\rm TL}}$ is the linear basis of the classical Temperley-Lieb algebra as computed by Jones in \cite{jo} and $\mathcal{E}_{d,n}(w)$ is a subset of $\{ 0 , \ldots , d-1 \}^n$ that describes the exponents of the $t_i$'s that correspond to the fixed braid word $w \in {\rm YTL}_{d,n}(u)$. For an explicit description of the set $\mathcal{E}_{d,n}(w)$, the reader is encouraged to consider \cite[Propositions~9 and 11]{ChPou}.
Subsequently, the dimension of the Yokonuma-Temperley-Lieb algebra can be computed, which is equal to:
\[
 {\rm dim}( {\rm YTL}_{d,n}(u)) = d c_n + \frac{d(d-1)}{2} \sum_{k=1}^{n-1} 	\dbinom{n}{k}^2, 
 \]
where $c_n$ is the $n^{th}$ Catalan number \cite[Proposition 4]{ChPou}.

By standard results in representation theory we have that the irreducible representations of ${\rm YTL}_{d,n}(u)$ are in bijection with those irreducible representations of ${\rm Y}_{d,n}(u)$ that respect the defining relation of ${\rm YTL}_{d,n}(u)$, which is $ g_{1,2} =0$. Specifically, the irreducible representations of ${\rm YTL}_{d,n}(u)$ are those representations of ${\rm Y}_{d,n}(u)$ who have at most two columns in total in the Young diagram of the parametrizing $d$-partition of $n$ \cite[Theorem~1]{ChPou}. In the following example, the first 3-partition of 5 parametrizes an irreducible representation of ${\rm YTL}_{3,5}(u)$ while the second one does not correspond to an irreducible representation of ${\rm YTL}_{3,5}(u)$:
\[
i. \left (\  \yng(1,1,1) \ ,  \yng(1,1) \  , \emptyset \ \right ) \qquad ii. \left ( \ \yng(3) \ ,   \yng(1,1) \ , \emptyset \right ).
\]

As mentioned in the introduction, the motivation behind the definition of a Temperley-Lieb type quotient from the Yokonuma-Hecke algebra was the construction of polynomial invariants for framed knots and links via the use of the trace ${\rm tr}_d$ of ${\rm Y}_{d,n}(u)$. Thus, one of the biggest challenges regarding the study of the algebra ${\rm YTL}_{d,n}(u)$ was the determination of the necessary and sufficient conditions for ${\rm tr}_d$ to factor through to the quotient algebra. By employing the methods that P. G\'erardin used to describe the full set of solutions of the ${\rm E}$-system \cite[Appendix]{jula}, the author together with Juyumaya, Kontogeorgis and Lambropoulou proved that following:
  \begin{theorem}[{\cite[Theorem 6]{gojukola}}]\label{thmy}
The trace ${\rm tr}_d$ passes to the quotient algebra ${\rm YTL}_{d,n}(u)$ if and only if the $x_i$'s  are solutions of the ${\rm E}$-system and one of the two cases holds: 
\begin{enumerate}
\item [(i)] 
the $x_\ell$'s are  $d^{th}$ roots of unity and  $z=-\frac{1}{u+1}$ or $z=-1$,
\item [(ii)] 
the $x_\ell$'s are the solutions of the {\rm E}-system that are parametrized by the set $D=\{ m_1, m_2 \, | \, 0 \leq m_1,m_2 \leq d-1\, \mbox{and } m_1\neq m_2\}$ and they are expressed as: 
\[
x_\ell= \frac{1}{2} \left (\chi_{m_1}(t^\ell) + \chi_{m_2}(t^\ell) \right ),  \quad 0 \leq \ell \leq d-1,
\]
where the  $\chi_k$'s denote the characters of the group $\mathbb{Z}/d\mathbb{Z}$. In this case we have that $z=-\frac{1}{2}$.
\end{enumerate}
\end{theorem}

Note that in both cases the $x_i$'s are solutions of the ${\rm E}$-system, as required by \cite{jula}, in order to proceed with defining link invariants. We do not take into consideration case $(i)$ for $z=-1$ and case $(ii)$, where $z= - \frac{1}{2}$, since crucial braiding information is lost and therefore they are of no topological interest \cite{gojukola}. The only remaining case of interest is case $(i)$ of Theorem~\ref{thmy}, where the $x_\ell$'s are the $d^{th}$ roots of unity and $z= - \frac{1}{u+1}$.  This implies that $E=1$ and $w=u$ in \eqref{gammainv}. So, by \cite{ChLa} and \cite{jo}, the invariant $\Delta_{d,s}(u,u)$  coincides with the Jones polynomial. For this reason, the algebra ${\rm YTL}_{d,n}(u)$ is discarded as a potential framization of the Temperley-Lieb algebra.

\subsection{{\it The Complex Reflection Temperley-Lieb algebra}} We move on now with presenting the second natural definition of a potential framization of the Temperley-Lieb algebra. For $n\geq 3$, we define the {\it Complex Reflection Temperley-Lieb} algebra, denoted by ${\rm CTL}_{d,n}(u)$, as the quotient of the algebra ${\rm Y}_{d,n}(u)$ over the ideal that is generated by the elements 
\begin{equation}\label{c12exp}
c_{i,i+1}:=  \sum_{\alpha, \beta, \gamma \in \mathbb{Z}/d\mathbb{Z}} t_i^\alpha t_{i+1}^\beta t_{i+2}^\gamma \ g_{i,i+1}.
\end{equation}

In analogy to the algebra ${\rm YTL}_{d,n}(u)$, the defining ideal of ${\rm CTL}_{d,n}(u)$ can be shown to be principal and is generated by the single element $c_{1,2}$. Further, for $d=1$, the algebra ${\rm CTL}_{1,n}(u)$ coincides with the algebra ${\rm TL}_n(u)$. The denomination Complex Reflection Temperley-Lieb algebra has to do with the fact that the underlying group of ${\rm CTL}_{d,n}(u)$ is isomorphic to the complex reflection group $G(d,1,3)$.

The Complex Reflection Temperley-Lieb algebra is isomorphic to a direct sum of matrix algebras over tensor products of Temperley-Lieb and Iwahori-Hecke algebras \cite[Theorem~5.8]{ChPou2}. This isomorphism, which we will denote by $\phi_n$, will lead to the determination of a linear basis for ${\rm CTL}_{d,n}(u)$. More precisely, there exists an explicit isomorphism:
\[
\phi_n \ : \ \bigoplus_{\mu \in {\rm Comp}_d(n)} {\rm Mat}_{m_\mu}  \big( {\rm TL}_{\mu_1}(u) \otimes {\rm H}_{\mu_2}(u) \otimes \ldots  \otimes {\rm H}_{\mu_d}(u) \big) \longrightarrow {\rm CTL}_{d,n}(u).
\]
Then the following set is a linear basis for ${\rm CTL}_{d,n}(u)$ \cite[Proposition~5.9]{ChPou2}:
{\small
\[
\left \{ \phi_n\left(b_1 b_2 \ldots b_d\ {\rm M}_{k,l} \right) \ \vert \  b_1 \in \mathcal{B}_{{\rm TL}_{\mu_1}(u)}, b_i \in  \mathcal{B}_{{\rm H}_{\mu_i}(u)} \mbox{ for all } i= 2,...,d,1\leq k,l \leq m_\mu, \mu \in {\rm Comp}_d(n) \right\},
\]
}
where $\mathcal{B}_{{\rm TL}_{\mu_1(u)}}$ is the linear basis of ${\rm TL}_{\mu_1(u)}$, $\mathcal{B}_{{\rm H}_{\mu_i}(u)}$ is the linear basis of ${\rm H}_{\mu_i}$,  ${\rm M}_{k,l}$ is the elementary $m_\mu \times m_\mu$ matrix with 1 in position$(k,l)$ and  $\mu \in {\rm Comp}_d (n)$ is a $d$-composition of $n$, that is, $\mu= (\mu_1,\, \mu_2,\,  \ldots , \mu_d ) \in \mathbb{N}^{d}$ such that $\mu_1 + \mu_2 + \ldots + \mu_d = n$. Counting the elements of the above basis one can derive the dimension of the algebra ${\rm CTL}_{d,n}(u)$ \cite[Theorem~5.5]{ChPou2}. Indeed, if $c_k:= \frac{1}{k+1} {2k \choose k}$ is the $k$-th Catalan number, we have that:
\[
{\rm dim}_{\mathbb{C}(u)}{\rm CTL}_{d,n}(u) = \sum_{\mu\, \in \, {\rm Comp}_d(n)} \left( \frac{n !}{\mu_1 ! \, \mu_2 ! \, \ldots \, \mu_d !} \right)^2 c_{\mu_1}\, \mu_2 ! \, \ldots \, \mu_d! 
\]

Let now $\lambda = (\lambda^{(1)} , \ldots , \lambda^{(d)})$ a $d$-partition of $n$. The irreducible representations of ${\rm CTL}_{d,n}(u)$ are those irreducible representations of ${\rm Y}_{d,n}(u)$ whose Young diagram of $\lambda^{(1)}$ has at most two columns \cite[Theorem~5.3]{ChPou2}. For instance, in the example given below, the first 2-Young diagram corresponds to an irreducible  representation of ${\rm CTL}_{2,9}(u)$ while the second one does not:
\[
i. \left ( \ \yng(2,2,1) \ ,  \ \yng(3,1) \ \right ) \qquad ii. \left ( \ \yng(3,1,1) \ ,  \ \yng(3,1) \ \right ).
\]

Next we present the necessary and sufficient conditions for the trace ${\rm tr}_d$ to factor through to the quotient algebra ${\rm CTL}_{d,n}(u)$. We have the following:
\begin{theorem}{{\cite[Theorem~7]{gojukola2}}}\label{ctlthm}
The trace ${\rm tr}_d$ passes to the quotient algebra ${\rm CTL}_{d,n}(u)$ if and only if the parameter $z$ and the $x_i$'s are related through the equation:
 \begin{equation}\label{ctlbasic}
 (u+1)z^2\sum_{k\in \mathbb{Z}/d\mathbb{Z}} x_k  +  (u+2)z\sum_{k\in \mathbb{Z}/d\mathbb{Z}} E^{(k)}+ \sum_{k \in \mathbb{Z}/d\mathbb{Z}}{\rm tr}(e_1^{(k)} e_2)=0.
 \end{equation}
\end{theorem}

Notice now that the conditions of Theorem~\ref{ctlthm} do not include any solutions of the ${\rm E}$-system. Thus, in order to obtain any well defined invariant from the algebras ${\rm CTL}_{d,n}(u)$ one has to impose the ${\rm E}$-condition on the trace parameters $x_i$. Even by doing so, ${\rm CTL}_{d,n}(u)$ does not deliver any new invariants for framed or classical oriented knots and links. We have the following:
\begin{proposition}[{\cite[Proposition~10]{gojukola2}}]
Let $X_D$ be a solution of the ${\rm E}$-system parametrized by the subset $D$ of $\mathbb{Z}/d\mathbb{Z}$. The invariants derived from the algebra ${\rm CTL}_{d,n}(u)$:
\begin{enumerate}
\item if $0 \in D$, they coincide with the invariants derived from the algebra ${\rm FTL}_{d,n}(u)$, 
\item if $0 \notin D$, they coincide with the invariants derived from the algebra ${\rm Y}_{d,n}(u)$.
\end{enumerate}
\end{proposition}

The above constitute the reasons for which the Complex Reflection Temperley-Lieb algebra is discarded as a potential candidate for the framization of the Temperley-Lieb algebra.

\subsection{{\it The Framization of the Temperley-Lieb algebra}}\label{ftldefsec}
For $n \geq 3$, the {\it Framization of the Temperley-Lieb algebra}, denoted by ${\rm FTL}_{d,n}(u)$, is defined as the quotient ${\rm Y}_{d,n}(u)$ over the two-sided ideal that is generated by the elements
\begin{equation}\label{ftlideal}
r_{i,i+1}:= \sum_{\alpha+\beta+\gamma=0} t_i^\alpha t_{i+1}^{\beta}  t_{i+2}^\gamma \ g_{i,i+1}.
\end{equation}
In analogy to the case of the other two quotient algebras, for $d=1$ the algebra ${\rm FTL}_{1,n}(u)$ coincides with ${\rm TL}_n(u)$. Additionally, the defining ideal of ${\rm FTL}_{d,n}(u)$ is principal and is generated by the element $r_{1,2}$.
Thus, in terms of generators and relations, ${\rm FTL}_{d,n}(u)$ is the $\mathbb{C}(u)$-algebra  generated by the set $\{t_1,\ldots , t_n ,$  $g_1 , \ldots , g_{n-1}\}$ whose elements are subject to the defining relations of ${\rm Y}_{d,n}(u)$ and the relation $r_{1,2}=0$.

As in the case of ${\rm CTL}_{d,n}(u)$, the determination of a linear basis for the Framization of the Temperley-Lieb algebra will emerge from an isomorphism theorem for ${\rm FTL}_{d,n}(u)$. More precisely, the quotient algebra ${\rm FTL}_{d,n}(u)$ is isomorphic to a direct sum of matrix algebras over tensor products of Temperley-Lieb algebras \cite[Theorem~4.3]{ChPou2}. There exists an explicit isomorphism of $\mathbb{C}(u)$-algebras:
\[
\widetilde{\phi}_n \ : \ \bigoplus_{\mu \in {\rm Comp}_d(n)} {\rm Mat}_{m_\mu}  \big( {\rm TL}_{\mu_1}(u) \otimes \ldots \otimes  {\rm TL}_{\mu_d}(u) \big ) \longrightarrow {\rm FTL}_{d,n}(u),
\]
then the following set is a linear basis for the algebra ${\rm FTL}_{d,n}(u)$:
\[
\left \{\widetilde{\phi}_n (b_1 \ldots b_d \  M_{k,l}) \ \vert \   b_i \in \mathcal{B}_{{\rm TL}_{\mu_i}(q)} \mbox{ for all } i = 1, \ldots d, 1 \leq k, l \leq m_{\mu}, \mu \in {\rm Comp}_d(n) \right \}.
\]
By using a counting argument one can derive the dimension of the algebra ${\rm FTL}_{d,n}(u)$, which is equal to \cite[Theorem~3.11]{ChPou2}:
\begin{equation}\label{dimftl}
{\rm dim}_{\mathbb{C}(u)}{\rm FTL}_{d,n}(u) = \sum_{\mu \, \in \, {\rm Comp}_d (n)} \left ( \frac{n!}{\mu_1! \, \mu_2 ! \, \ldots \mu_d !} \right )^2 c_{\mu_1}\, c_{\mu_2}\ \ldots \, c_{\mu_d}.  
\end{equation} 
The irreducible representations of ${\rm FTL}_{d,n}(u)$ are those irreducible representations of ${\rm Y}_{d,n}(u)$ whose Young diagram of $\lambda^{(i)}$ has at most two columns, for $i=1, 2 ,\ldots , d$. As in the previous examples, the first of the following 3-Young diagrams describes an irreducible representation of ${\rm FTL}_{3,7}(u)$ while the second does not:
\[
i. \left ( \ \yng(2,1)\ , \ \yng(2,1)\ ,\ \yng(1) \ \right ) \qquad ii. \left ( \ \yng(2,1)\ , \ \yng(3)\ ,\ \yng(1) \ \right ).
\]
We move on now to the necessary and sufficient conditions so that ${\rm tr}_d$ factors through to ${\rm FTL}_{d,n}(u)$.
\begin{theorem}[{\cite[Theorem~6]{gojukola2}}]\label{ftlthm}
The trace {\rm tr } passes to ${\rm FTL}_{d,n}(u)$ if and only if the parameters of the trace {\rm tr} satisfy: 
\[
x_k = -z \left(\sum_{m\in {\rm Sup}_1}\chi_{ m}(t^{k}) + (u+1)\sum_{m\in {\rm Sup}_2}\chi_{ m}(t^{k}) \right)
\quad \text{and}\quad  
z=-\frac{1}{\vert {\rm Sup_1}\vert + (u+1)\vert {\rm Sup_2}\vert  },
\]
where $\chi_m$ are the characters of the group $\mathbb{Z}/d\mathbb{Z}$, ${\rm Sup}_1\sqcup \rm{Sup}_2$ (disjoint union) is the support of the Fourier transform of $x$, and $x$ is the complex function on $\mathbb{Z}/d\mathbb{Z}$,   
that maps $0$ to $1$ and $k$ to the trace parameter  $x_k$.

\end{theorem}
The intrinsic difference with the other two quotient algebras lies in the fact that the necessary and sufficient conditions of Theorem~\ref{ftlthm} include all solutions of the ${\rm E}$-system. This observation is the main reason that led to the consideration of the quotient algebra ${\rm FTL}_{d,n}(u)$ as the most natural non-trivial analogue of the Temperley-Lieb algebra in the context of framization of knot algebras. If one lets either ${\rm Sup}_1$ or ${\rm Sup}_2$ to be the empty set, then the trace parameters $x_k$ comprise a solution of the ${\rm E}$-system. In this context, if ${\rm Sup}_1$ is the empty set  then $z= - \frac{1}{(u+1)|{\rm Sup}_2|}$ while if ${\rm Sup}_2$ is the empty set then $z=-1/|{\rm Sup}_1|$ \cite[Corollary~3]{gojukola2}. Since for defining invariants for oriented (framed) knots and links only the cardinal $|D|$ of the parametrizing set $D$ of a solution is needed, the solutions mentioned above cover all the possibilities. 
We do not take into consideration the case where ${\rm Sup}_2 = \emptyset$ and  $z=-1/|{\rm Sup}_1|$ since important topological information is lost and thus basic pairs of knots are not distinguished \cite[Remark~7]{gojukola2}. For the remaining case, let $X_D$ be a solution of the ${\rm E}$-system, parametrized by the non-empty subset $D = {\rm Sup}_2$ of $\mathbb{Z}/d\mathbb{Z}$ and let $z = -\frac{1}{(u+1)|D|}$. We obtain from $\Gamma_{d,D}(w,u)$ the following new 1-variable framed link invariants:
\begin{equation}
\begin{array}{crcl}
 \Gamma_{d,D}(u, u)(\widehat{\alpha})&:=& \left (- \frac{(1+ u)|D|}{\sqrt{u} } \right)^{n-1} \left(\sqrt{u}\right)^{\varepsilon(\alpha)} {\rm tr}_{d,D}\left(\gamma(\alpha)\right) \label{frjon},
\end{array}
\end{equation}
for any $\alpha \in \cup_{\infty}\mathcal{F}_n$. Further, in analogy to  the invariants of $\Gamma_{d,D}(w,u)$, if we restrict to framed links with all framings zero, we obtain from  $\Gamma_{d,D}(u,u)$ new 1-variable invariants of classical links $\Delta_{d,D}(u,u)$. Additionally, for $d=1$ the invariant $\Gamma_{d,D}(u,u)$ coincides with the Jones polynomial. 

\section{Comparisons and generalizations}\label{compgen}
In this section we will present the comparisons of the invariants $\Theta_d$ and $\theta_d$ to the Homflypt and the Jones polynomials respectively, and we will give generalizations for both of them.
\subsection{{\it The invariants $\Theta_d$ and their generalization}} In a recent development \cite{ChJuKaLa} it was proved that the classical link invariants derived from the Yokonuma-Hecke algebra are {\it not topologically equivalent to the Homflypt polynomial on links} while they are topologically equivalent to the Homflypt on knots. This was achieved by considering a different presentation for the algebra ${\rm Y}_{d,n}$ with parameter $q$ instead of $u$ and a different quadratic relation. More precisely, the algebra ${\rm Y}_{d,n}(q)$ is defined as the $\mathbb{C}(q)$-algebra that is generated by the elements $g_1^\prime , \ldots , g_{n-1}^\prime , t_1 ,\ldots , t_n$, which satify all relations of ${\rm Y}_{d,n}(u)$ except for the quadratic relation that is replaced with the following:
\begin{equation}\label{newquad}
(g^\prime_i)^{\,2} = 1+ (q-q^{-1}) e_i g^\prime_i.
\end{equation}

One can obtain this presentation from the one given in Section~\ref{yhsec} by taking $u=q^2$ and 
\[
g_i=g^\prime_i+(q-1)e_ig^\prime_i  \quad \mbox{(or, equivalently, }\, g^\prime_i=g_i+(q^{-1}-1)e_i g_i). 
\]

Thus, the following invariants of classical links were derived \cite{ChJuKaLa}:
\begin{equation}\label{Thetadinv}
\Theta_d(q,\lambda_d)(\widehat{\alpha}) = \left( \frac{1-\lambda_d}{\sqrt{\lambda_d}(q-q^{-1})E_D} \right)^{n-1}\sqrt{\lambda_d}^{\varepsilon(\alpha)}{\rm tr}_{d,D}(\delta(\alpha)),
\end{equation}
where $\alpha \in \cup_\infty B_n$,  $E_D=1/d$, $\varepsilon(a)$ is as in \eqref{gammainv}, $\delta$ is the natural homomorphism $\mathbb{C}(q)B_n \rightarrow {\rm Y}_{d,n}(q)$ and $\lambda_d=\frac{{z'} -(q-q^{-1})E_D}{z'}$ is the re-scaling factor for the trace ${\rm tr}_{d,D}$.  

The invariants $\Theta_d$ depend only on $d\in \mathbb{N}$, that is, the cardinal of the subset $D$ that parametrizes the solution of the ${\rm E}$-system \cite[Proposition~4.6]{ChJuKaLa}. Furthermore, the choice of the new presentation for ${\rm Y}_{d,n}$ revealed that the invariants $\Theta_d$ satisfy the Homflypt skein relation on crossings between different components of a link $L$ \cite[Proposition~6.8]{ChJuKaLa}. Using this, one can prove that the invariants $\Theta_d$ distinguish more pairs of Homflypt equivalent pair of non-isotopic oriented classical links \cite[Section~7.2]{ChJuKaLa} and thus that $\Theta_d$ are {\it not topologically equivalent to the Homflypt polynomial on links} \cite[Theorem~7.3]{ChJuKaLa}.

In \cite{ChJuKaLa} it has been shown skein-theoretically that the invariants for classical links $\Theta_d$ generalize to a new 3-variable invariant $\Theta (q, \lambda, E)$ for classical oriented links that can be defined uniquely by the following two rules:
 \begin{enumerate}
\item On crossings between different components of an oriented classical link $L$ the skein relation of the Homflypt polynomial holds:
\[
\frac{1}{\sqrt{\lambda_D}} \, \Theta(L_+) - \sqrt{\lambda_D}\, \Theta (L_{-}) = (q-q^{-1})\,  \Theta(L_0),
\]
where $L_+$, $L_-$ and $L_0$ is a Conway triple. 
\item  For a disjoint union of $\mathcal{K} = \sqcup_{i=1}^r K_i$ of $r$ knots, with $r>1$, it holds that:
\[
\Theta (\mathcal{K})= E^{1-r} \prod_{i=1}^{r} P(K_i),
\]
where $P(K_i)$ is the value of the Homflypt polynomial on $K_i$.
\end{enumerate}
Algebraically, the well-definedness of the invariant $\Theta$ can be proved by using the {\it the algebra of braids and ties},  $\mathcal{E}_n(q)$ \cite{AiJu}. The algebra $\mathcal{E}_n(q)$ supports a unique Markov trace $\rho$ that gives rise to a 3-variable invariant for tied links $\overline \Theta (q, \lambda, E)$ which, in turn, restricts to an invariant of classical oriented links $ \Theta (q, \lambda, E)$ \cite{AiJu1, AiJu2}.   Alternatively, one can use the fact that $\mathcal{E}_n(q)$  is isomorphic to the subalgebra ${\rm Y}_{d,n}^{{(\rm br)}}(q)$ of ${\rm Y}_{d,n}(q)$ that is generated only by the $g_i$'s \cite{ERH}.  Note now that when computing the specialized trace ${\rm tr}_{d,D}$ of a braid word in $B_n$, the framing generators appear only when applying the quadratic or the inverse relation and only in the form of the idempotents $e_i$. 
In this case and by the ${\rm E}$-condition, the last rule of  the specialized trace: ${\rm tr}_{d,D}(at_{n+1}^s) =  {\rm x}_s {\rm tr}_{d,D}(a)$, for $s = 1, \ldots , d-1$,  can be substituted by the following two rules \cite[Theorem~4.3]{ChJuKaLa}:
 \[
 {\rm tr}_{d,D}(a e_n) = E_D \, {\rm tr}_{d,D}(a) \quad \mbox{and} \quad {\rm tr}_{d,D}(a e_n g_n) = z\, {\rm tr}_{d,D}(a),
 \]
where $D$ is the non-empty subset of $\mathbb{Z}/d\mathbb{Z}$ that parametrizes a solution of the ${\rm E}$-system. Consequently, if  $E_D$ is considered as an indeterminate, the specialized trace ${\rm tr}_{d,D}$ on ${\rm Y}^{(\rm br)}_{d,n}(q)$ is well-defined since it coincides with the trace $\rho$ on $\mathcal{E}_n(q)$ and, therefore, the invariant $\Theta$ can be constructed directly through ${\rm Y}_{d,n}^{(\rm br)}(q)$ \cite[Remark~4.18]{ChJuKaLa}. Conversely, one can recover the invariants $\Theta_d$ from $\Theta$ by specializing $E=1/d$, $d \in \mathbb{N}$.

   A self-contained diagrammatic proof for the well-definedness of the invariant $\Theta$ has been given in \cite{kaula}. The invariant $\Theta$ distinguishes more pairs of non isotopic oriented links than the Homflypt polynomial and thus it is stronger than the Homflypt. We note also that, $\Theta$ is not topologically equivalent to the Homflypt or the Kauffman polynomials. 
   
Finally, it is worth noting that the invariant $\Theta$ can be described by the following closed combinatorial formula, namely:
\begin{theorem}[{\cite[Appendix]{ChJuKaLa}}]\label{licktheta}
Let  $L$ be an oriented link with $n$ components, then:

\begin{equation}\label{clform}
\Theta(q,\lambda,E)(L) = \sum_{k=1}^m \mu^{k-1} E_k \sum_{\pi} \lambda^{\nu(\pi)} P(\pi L), 
\end{equation}
where the second summation is over all partitions of $\pi$ of the components of $L$ into $k$ (unordered) subsets and $P(\pi L)$ denotes the product of the Homflypt polynomial of the $k$ sublinks of $L$ defined by $\pi$. Furthermore, $\nu(\pi)$ is the sum of all linking numbers of pairs of components of $L$ that are distinct sets of $\pi$, $E_k = (E^{-1} -1 ) (E^{-1} -2 ) \ldots (E^{-1} - k+1)$, with $E_1=1$ and $\mu = \frac{\lambda^{-1/2} - \lambda^{1/2}}{q-q^{-1}}$
\end{theorem}

\subsection{{\it The invariants $\theta_d$ and their generalization}}  By adjusting the algebra ${\rm FTL}_{d,n}$ to the presentation that has parameter $q$ and involves the quadratic relation \eqref{newquad}, one can compare the derived invariants for classical oriented links to the Jones polynomial. In this context, the generator of the defining principal ideal of ${\rm FTL}_{d,n}$ is transformed to the following element of ${\rm Y}_{d,n}(q)$:
 \[
  e_1 e_2 \Big( 1 + q(g^\prime_1 + g^\prime_2) + q^2 (g^\prime_1 g^\prime_2 + g^\prime_2 g^\prime_1) + q^3 g^\prime_1 g^\prime_2 g^\prime_1 \Big).
 \]
Note that the ${\rm E}$-system and its solutions remain unaffected by this change of presentations. The values for the trace parameters $z$, however, are transformed to the following:
\[
z^\prime=-\frac{q^{-1}E_D}{q^2+1} \quad \mbox{or} \quad z^\prime=-q^{-1} E_D.
\]
The parameters  $z$ and $z^\prime$ are related through the equation: $z = q z^\prime$. Again, the value $z^\prime=-q^{-1} E_D$ is discarded. For the remaining values for $z^\prime$, we obtain from \eqref{Thetadinv} the following 1-variable specialization of $\Theta_d$:
\[
\theta_d(q)(\widehat{\alpha}) := \left ( -\frac{1+q^2}{qE_D} \right)^{n-1} q^{2\varepsilon(\alpha)} {\rm tr}_{d,D} (\delta(a)) = \Theta_d(q, q^4)(\widehat{\alpha}),
\]
where $\alpha \in \cup_\infty B_n$, $d$ and $E_D$, $\varepsilon(a)$ and $\delta$ are as in \eqref{Thetadinv}. The invariants $\theta_d$ were proven to be topologically equivalent to the Jones polynomial on knots \cite[Proposition~11]{gojukola2}, however, they are {\it topologically not equivalent to the Jones polynomial on links} \cite[Theorem~9]{gojukola2}.

In \cite{gola1} the author together with S. Lambropoulou has shown that the invariants $\theta_d$ generalize to a new 2-variable invariant $\theta$ for classical links. This generalization can be proved either algebraically or diagrammatically. Algebraically, this can be shown in two different ways. The first way is to consider the {\it partition Temperley-Lieb} algebra, ${\rm PTL}_n(q)$, which is a quotient of $\mathcal{E}_n(q)$ and determine the necessary and sufficient conditions such that the uinque Markov trace $\rho$ of $\mathcal{E}_n(q)$ factors through to ${\rm PTL}_n(q)$. These conditions give rise to a 2-variable invariant for classical links, $\theta(q,E)$ \cite[Definition~1]{gola1}, that for $E=1/d$ coincides with $\theta_d$.
Alternatively, one can  show that, for $d \geq n$, the subalgebra ${\rm FTL}_{d,n}^{{\rm (br)}}(q)$ of ${\rm FTL}_{d,n}$ that is generated only by the braiding generators $g_i$ is isomorphic to ${\rm PTL}_n(q)$ \cite[Proposition~5]{gola}.  Diagrammatically, one may consider the skein-theoretic definition of $\Theta(q,\lambda,E)$ and specialize $\lambda=q^4$. Thus, we obtain the following:

\begin{theorem}[{\cite[Theorem~6]{gola1}}]\label{mainthm}
Let $q, E$ be indeterminates. There exists a unique ambient isotopy invariant of classical oriented links
\[
\theta : \mathcal{L} \rightarrow \mathbb{C}[q^{\pm 1} , E^{\pm 1}]
\]
defined by the following rules:
\begin{enumerate}
\item On crossings involving different components the following skein relation holds:
\[
q^{-2}\, \theta (L_+) - q^2\, \theta (L_-) = (q - q^{-1})\, \theta (L_0),
\]
where $L_+$, $L_-$ and $L_0$ constitute a Conway triple.
\item For a disjoint union $\mathcal{K} = \sqcup_{i=1}^r K_i$ of $r$ knots, with $r>1$, it holds that:
\[
\theta (\mathcal{K}) = E^{1-r} \prod_{i=1}^{r}V(K_i),
\]
where $V(K_i)$ is the value of the Jones polynomial on $K_i$.
\end{enumerate}
\end{theorem}

All the properties of the invariant $\Theta$ carry through to $\theta$ \cite{gola1} and so the invariant $\theta$ distinguishes the same pairs of non-isotopic oriented classical links as $\Theta$. More precisely, in \cite{ChJuKaLa} six pairs of Homflypt-equivalent non-isotopic oriented classical links were found to be distinguished by the invariants $\Theta(q,\lambda,E)$, which are all still distinguished by $\theta$. Indeed we have that:

{\small\begin{align*}
&\theta(L11n358\{0,1\})-\theta(L11n418\{0,0\}) = \frac{(1-E) (q-1)^5 (q+1)^5 (q^2+1)(q^2+q+1)(q^2-q+1)}{E\, q^{18}}
\\
&\theta(L11a467\{0,1\})-\theta(L11a527\{0,0\}) =  \frac{(1-E) (q-1)^5 (q+1)^5 (q^2+1)(q^2+q+1)(q^2-q+1)}{E\, q^{18}}
\\
&\theta(L11n325\{1,1\})-\theta(L11n424\{0,0\}) = \frac{(E-1) (q-1)^5 (q+1)^5 (q^2+1)(q^2+q+1)(q^2-q+1)}{E\, q^{14}}\\
&\theta(L10n79\{1,1\})-\theta(L10n95\{1,0\}) = \frac{(E-1) (q^2-1)^3 (q^8+2\,q^6+2\,q^4-1)}{E \, q^{18}}\\
&\theta(L11a404\{1,1\})-\theta(L11a428\{0,1\}) = \frac{(1-E) (q-1)^3(q+1)^3(q^2+1)(q^4+1)(q^6-q^4+1)}{E\, q^4}
\\  
&\theta(L10n76\{1,1\})-\theta(L11n425\{1,0\}) =  \frac{(E-1) (q-1)^3(q+1)^3(q^2+1)(q^4+1)}{E\, q^{10}}.\\
\end{align*}}

The invariant $\theta(q,E)$ is not topologically equivalent to the Homflypt or the Kauffman polynomials, it includes the family of invariants $\{ \theta_d \}_{d\in\mathbb{N}}$ as well as the Jones polynomial and hence it is {\it stronger than the Jones polynomial} \cite[Theorem~7]{gola1}.

Finally, the invariant $\theta$ can be described by a closed combinatorial formula, which is a corollary of Theorem~\ref{licktheta}. Indeed we have:

\begin{corollary}\label{corol}
Let  $L$ be an oriented link with $n$ components. Then:
\[
\theta(q,E)(L) = \sum_{k=1}^m (-1)^{k-1} (q+q^{-1})^{k-1} E_k \sum_{\pi} \lambda^{\nu(\pi)} V(\pi L), 
\]
where $\pi$, $\nu(\pi)$, and $E_k$ are as in Theorem~\ref{licktheta}, and $V(\pi L)$ denotes the product of the Jones polynomial of the $k$ sublinks of $L$ defined by $\pi$. 
\end{corollary}

From Corollary~\ref{corol} it is clear that the invariant $\theta$ depends on the orientations of the components of the link $L$, thus making it impossible to relate $\theta$ to the Kauffman bracket polynomial. However, as shown in \cite[Theorem~7]{gola1}, $\theta$ can be expressed in terms of the oriented extension of the bracket polynomial. In particular, the author together with S. Lambropoulou defined in \cite{gola1} the ambient isotopy link invariant $\left \{ \left \{ L \right \}\right \}$ of the  link diagram $L$ by the following two rules: \\
\smallbreak
\noindent (1) For a disjoint union $\mathcal{K}^r : = \sqcup_{i=1}^r K_i$, of $r$ knots with $r \geq 1$, we have that:
\begin{equation}\label{init}
\left \{ \left \{ \mathcal{K}^r\right \}\right \}:= E^{1-r} \prod_{i=1}^r V (K_i),
\end{equation}
(2) On crossings involving different components the skein relation of the Jones polynomial holds, namely:
\begin{equation}\label{curlbrack}
q^{-2} \, \left \{  \left \{ \raisebox{-.1cm}{\begin{tikzpicture}[scale=.2]
\draw [line width=0.35mm, draw=purple]  (-1,-1)-- (-0.22,-0.22);
\draw  [line width=0.35mm, draw=blue](-1,1)--(0,0);
\draw  [line width=0.35mm, draw=purple] (0.22,0.22) -- (1,1)[->];
\draw [line width=0.35mm, draw=blue]   (0,0) -- +(1,-1)[->];
\end{tikzpicture}}\, \right \} \right \} - q^2 \, 
 \big \{  \big \{ \raisebox{-.1cm}{\begin{tikzpicture}[scale=.2]
\draw  [line width=0.35mm, draw=purple] (-1,-1)-- (0,0) ;
\draw [line width=0.35mm, draw=blue] (-1,1)--(-0.22,0.22);
\draw [line width=0.35mm, draw=purple] (0,0) -- (1,1)[->];
\draw [line width=0.35mm, draw=blue]   (0.22,-0.22) -- +(.8,-.8)[->];
\end{tikzpicture}} \,  \big\} \big\}
= (q- q^{-1} ) \,
\left \{\left \{ \, \raisebox{-.1cm}{\begin{tikzpicture}[scale=.2, mydeco/.style = {decoration = {markings, 
                                                       mark = at position #1 with {\arrow{>}}}
                                       }]
\draw [purple, line width=0.35mm, postaction = {mydeco=.6 ,decorate}] plot [smooth, tension=2] coordinates { (-1,.8) (0, 0.5) (1,.8)};
\draw [purple, line width=0.35mm, postaction = {mydeco=.6 ,decorate}] plot [smooth, tension=2] coordinates { (-1,-.8) (0, -0.5) (1,-.8)};
\end{tikzpicture}}\, \right \}\right \}.
\end{equation}
Comparing \eqref{init} and \eqref{curlbrack} to Theorem~\ref{mainthm}, we deduce that $\left \{\left \{L \right \}\right \}$ coincides with  the invariant $\theta(q,E)$.

\bibliography{bibliography.bib}{}

\begin{thebibliography}{10}

\bibitem{AiJu}
{\sc F.~Aicardi and J.~Juyumaya}, {\em An algebra involving braids and ties}.
\newblock Preprint ICTP IC/2000/179, Trieste, 2000.

\bibitem{AiJu1}
\leavevmode\vrule height 2pt depth -1.6pt width 23pt, {\em Markov trace on the
  algebra of braids and ties}, Moscow Math. J., 16 (2016), pp.~397--431.

\bibitem{AiJu2}
\leavevmode\vrule height 2pt depth -1.6pt width 23pt, {\em {Tied links}}.
\newblock to appear in Journal of Knot Theory and its Ramifications. {See also}
  arXiv:1503.00527, 2016.

\bibitem{ChPoulain}
{\sc M.~Chlouveraki and L.~P. D'Andecy}, {\em Representation theory of the
  {Yokonuma-Hecke}}, Advances in Mathematics, 259 (2014), pp.~134--172.

\bibitem{ChJuKaLa}
{\sc M.~Chlouveraki, J.~Juyumaya, K.~Karvounis, and S.~Lambropoulou}, {\em
  Identifying the invariants for classical knots and links from the
  {Yokonuma-Hecke} algebras}, submitted for publication. {See also}
  arXiv:1505.06666,  (2015).

\bibitem{ChLa}
{\sc M.~Chlouveraki and S.~Lambropoulou}, {\em {The {Yokonuma-Hecke} algebras
  and the {Homflypt} polynomial}}, {J. Knot Theory and Its Ramifications}, 22
  (2013).

\bibitem{ChPou}
{\sc M.~Chlouveraki and G.~Pouchin}, {\em Determination of the representations
  and a basis for the {Yokonuma-Temperley-Lieb} algebra}, Algebras and
  Representation Theory, 18 (2015).

\bibitem{ChPou2}
\leavevmode\vrule height 2pt depth -1.6pt width 23pt, {\em {Representation
  theory and an isomorphism theorem for the Framisation of the {Temperley-Lieb}
  algebra}}, to appear in Mathematische Zeitschrift. See also
  arXiv:1503.03396v2,  (2016).

\bibitem{ERH}
{\sc J.~Espinoza and S.~Ryom-Hansen}, {\em Cell structures for the
  {Yokonuma-Hecke} algebra and the algebra of braids and ties}, submitted for
  publication. See also arXiv:1506.00715.,  (2016).

\bibitem{homfly}
{\sc P.~Freyd, D.~Yetter, J.~Hoste, W.~Lickorish, K.~Millett, and A.~Ocneanu},
  {\em A new polynomial invariant of knots and links}, Bull. AMS, 12 (1985),
  pp.~239--246.

\bibitem{go}
{\sc D.~Goundaroulis}, {\em Framization of the Temperley-Lieb algebra and
  related link invariants}, PhD thesis, Department of Mathematics, National
  Technical University of Athens, 1 2014.

\bibitem{gojukola}
{\sc D.~Goundaroulis, J.~Juyumaya, A.~Kontogeorgis, and S.~Lambropoulou}, {\em
  {The Yokonuma-Temperley-Lieb Algebra}}, Banach Center Pub., 103 (2014),
  pp.~73--95.

\bibitem{gojukola2}
\leavevmode\vrule height 2pt depth -1.6pt width 23pt, {\em {Framization of the
  Temperley-Lieb Algebra}}, to appear in {Mathematical Research Letters}. See
  also arXiv:1304.7440v3,  (2016).

\bibitem{gola}
{\sc D.~Goundaroulis and S.~Lambropoulou}, {\em {Classical link invariants from
  the framizations of the Iwahori-Hecke algebra and the Temperley-Lieb algebra
  of type $A$}}, to appear in {Journal of Knot Theory and its Ramificiations}.
  See also {arXiv}:1602.07203,  (2016).

\bibitem{gola1}
\leavevmode\vrule height 2pt depth -1.6pt width 23pt, {\em A new two-variable
  generalization of the {Jones} polynomial}.
\newblock Submitted for publication, see also arXiv:1608.01812 [math.GT], 2016.

\bibitem{jo}
{\sc V.~Jones}, {\em Hecke algebra representations of braid groups and link
  polynomials}, Annals of Mathematics, 126 (1987), pp.~335--388.

\bibitem{jusur}
{\sc J.~Juyumaya}, {\em Sur les nouveaux g\'en\'erateurs de l'alg\`ebre de
  {Hecke} ${\mathcal h}(g,u,1)$}, J. Algebra, 204 (1998), pp.~40--68.

\bibitem{ju}
\leavevmode\vrule height 2pt depth -1.6pt width 23pt, {\em Markov trace on the
  {Yokonuma-Hecke} algebra}, J. Knot Theory and Its Ramifications, 13 (2004),
  pp.~25--39.

\bibitem{jukan}
{\sc J.~Juyumaya and S.~Kannan}, {\em Braid relations in the {Yokonuma-Hecke}
  algebra}, J. Algebra, 239 (2001), pp.~272--297.

\bibitem{jula2}
{\sc J.~Juyumaya and S.~Lambropoulou}, {\em $p$-adic framed braids}, Topology
  and its Applications, 154 (2007), pp.~1804--1826.

\bibitem{jula4}
\leavevmode\vrule height 2pt depth -1.6pt width 23pt, {\em An adelic extension
  of the jones polynomial}, in The mathematics of knots, M.~Banagl and
  D.~Vogel, eds., Contributions in the Mathematical and Computational Sciences,
  Vol. 1, Springer, 2009, pp.~825--840.

\bibitem{jula3}
\leavevmode\vrule height 2pt depth -1.6pt width 23pt, {\em An invariant for
  singular knots}, J. Knot Theory and Its Ramifications, 18 (2009),
  pp.~825--840.

\bibitem{jula5}
\leavevmode\vrule height 2pt depth -1.6pt width 23pt, {\em Modular framization
  of the {BMW} algebra}.
\newblock arXiv:1007.0092v1 [math.GT], 2013.

\bibitem{jula}
\leavevmode\vrule height 2pt depth -1.6pt width 23pt, {\em $p$-adic framed
  braids {II}}, Advances in Mathematics, 234 (2013), pp.~149--191.

\bibitem{jula6}
\leavevmode\vrule height 2pt depth -1.6pt width 23pt, {\em On the framization
  of knot algebras}, in {New Ideas in Low-dimensional Topology}, L.~Kauffman
  and V.~Manturov, eds., Series on Knots and everything, World Scientific,
  2014.

\bibitem{kaula}
{\sc L.~H. Kauffman and S.~Lambropoulou}, {\em New invariants of links and
  their state sum models}.
\newblock In preparation, private communication, 2016.

\bibitem{ks}
{\sc K.~Ko and L.~Smolinsky}, {\em The framed braid group and $3$-manifolds},
  Proceedings of the AMS, 115 (1992), pp.~541--551.

\bibitem{marin}
{\sc I.~Marin}, {\em Artin groups and the {Yokonuma-Hecke} algebra}.
\newblock arXiv:1601.03191, 2016.

\bibitem{pt}
{\sc J.~H. Przytycki and P.~Traczyk}, {\em {Invariants of links of Conway
  type}}, Kobe J. Math., 4 (1987), pp.~115--139.

\bibitem{yo}
{\sc T.~Yokonuma}, {\em {Sur la structure des anneux de Hecke d'un group de
  Chevalley fin}}, C.R. Acad. Sc. Paris, 264 (1967), pp.~344--347.

\end{thebibliography}
\bibliographystyle{siam}
\end{document}